\documentclass[11pt]{amsart}
\usepackage{amssymb, latexsym}
\theoremstyle{plain}
\newtheorem{theorem}{Theorem}

\newtheorem {lemma}{Lemma}
\newtheorem{proposition}{Proposition}

\newtheorem*{1'}{Theorem 1-Bessel}
\newtheorem*{P2'}{Proposition 2-Bessel}
\newtheorem*{P3'}{Proposition 3-Bessel}
\newtheorem*{P4'}{Proposition 4-Bessel}
\newtheorem*{C1'}{Corollary 1-Bessel}

\newtheorem*{2'}{Theorem 2-Bessel}
\newtheorem*{3'}{Theorem 3-Bessel}

\theoremstyle{remark}

\newtheorem*{Remark 1}{Remark 1}
\newtheorem*{Remark 2}{Remark 2}
\newtheorem*{Remark 3}{Remark 3}
\newtheorem*{Remark 4}{Remark 4}

\numberwithin{equation}{section}
\renewcommand{\baselinestretch}{1.4}

\begin{document}

\title[The secretary problem with pattern-avoiding permutations ] {The secretary problem with items arriving according to a random permutation avoiding a pattern of length three}

\author{Ross G. Pinsky and Tomer Zilca}

%\noindent  pinsky@math.technion.ac.il\ \ \ \ tel: 972-4-829-4083\ \ \  fax: 972-4-829-3388

\address{Department of Mathematics\\
Technion---Israel Institute of Technology\\
Haifa, 32000\\ Israel}
\email{ pinsky@technion.ac.il; tomerzilca@campus.techion.ac.il}

\urladdr{http://www.math.technion.ac.il/~pinsky/}

\date{}

\begin{abstract}
In the classical secretary problem, $n$ ranked items arrive one by one, and each item's rank relative to its predecessors is noted. The observer  must select or reject each item as it arrives, with the object of selecting the item of highest rank.
For $M_n\in\{0,1,\cdots, n-1\}$, let $\mathcal{S}(n,M_n)$ denote the  strategy whereby the observer rejects the first $M_n$ items, and then selects the first later-arriving item whose rank is higher than that of any
of the first $M_n$ items (if such an item exists).
If the ranked items arrive in a uniformly random order, it is well-known that the limiting optimal probability of success is $\frac1e$, which occurs if
$M_n\sim\frac ne$.
It has been shown that when   the ranked items arrive according to certain non-uniform distributions
on the set of permutations, $\frac1e$ serves as a lower bound for the optimal probability of success.
 There is a fundamental reason for this phenomenon, which we mention. In this note, we consider certain distributions on permutations for which that fundamental reason does not apply.
We begin by noting a cooked-up class of distributions for which $\mathcal{S}(n,M)$ yields the lowest possible probability of success---namely $\frac1n$, for all $M$.
We then consider a certain more natural class of distributions; namely  the uniform distribution over all permutations avoiding a particular pattern of length three.
In the case of the pattern 231 or  132,
 for \it any \rm\ choice of $M_n$, the  strategy $\mathcal{S}(n,M_n)$ yields the very same probability of success; namely $\frac{n+1}{2(2n-1)}$, which gives a limiting probability of success of $\frac14$.
For the pattern 213, the optimal strategy is obtained for  $M\in\{0,1\}$,  also yielding a limiting probability of success of $\frac14$.
On the other hand, for the pattern
123, the optimal strategy is obtained for $M=1$, yielding a limiting probability of success of $\frac34$.
For the other two patterns,
 312 and 321, an optimal strategy will yield a limiting probability of at least $\frac7{16}$.

\end{abstract}

\maketitle
\section{Introduction and Statement of Results}\label{intro}
\renewcommand{\baselinestretch}{1.3}

Recall the classical secretary problem: For $n\in\mathbb{N}$, a set of $n$ ranked items is revealed, one item at a time, to an observer whose objective is to select the item with
the highest rank. The order of the items is completely random; that is, each of the  $n!$ permutations of the ranks  $[n]=\{1,\cdots, n\}$ is equally likely.
At each stage, the observer only knows the relative ranks of the items that have arrived thus far, and must  either select the current item, in which case the process terminates, or reject it and continue to the next item. If the observer rejects the first $n-1$ items, then the $n$th and final item to arrive must be accepted.
As is very well known, asymptotically as $n\to\infty$, the optimal  strategy is to reject the first $M_n$ items, where $M_n\sim \frac ne$, and then to select the first
later-arriving item whose rank is higher than that of any of the first $M_n$ items (if such an item exists).
The limiting probability of successfully selecting the item of highest rank is $\frac1e$.

In \cite{P22} one of the authors considered the secretary problem when the random order in which the items arrive is not according to  the uniform distribution
on $S_n$, the set of  permutations of $[n]$, but rather according to a distribution from the family of   Mallows distributions, which are obtained from the uniform distribution by exponential tilting via the inversion statistic.  For certain ranges of the parameter defining the Mallows distribution, the limiting optimal probability of success using an optimal strategy  is $\frac1e$, as in the classical case, and
for other ranges it is greater than $\frac1e$.  This same phenomenon  also holds when the random order of arrival is according to a distribution obtained by exponential tilting via the left-to-right minimum statistic
\cite{Pman}.

 The fact that in the above-mentioned cases the optimal asymptotic  probability of success never dips below $\frac1e$ can be explained by a result of Bruss \cite{B00}.
For $n\in\mathbb{N}$, let $\{I_j\}_{j=1}^n$ be a sequence of independent indicator functions, which are observed sequentially. The observer's objective is to stop at the last $k$ for which $I_k=1$. Let $p_j\in[0,1]$ denote the probability that $I_j=1$.
One of the results of that paper is that an optimal strategy as $n\to\infty$   yields an optimal limiting probability of at least $\frac1e$, as long as
of $\sum_{j=1}^\infty \frac{p_j}{1-p_j}\ge 1$.

This result of Bruss can be applied to the classical secretary problem. Indeed, let $I_k$ be equal to 1 or 0 according to whether or not the $k$th item is the highest ranked item among the first $k$ items.
Since $p_1=1$, the summation condition above is automatically satisfied.
It is easy to check that
the $\{I_k\}_{k=1}^n$ are independent under the uniform distribution. It turns out that this independence also holds under the two aforementioned distributions. (The proof uses a certain online construction for each of these random permutations, along with the fact that if $\sigma\in S_n$ is such a random permutation, then $\sigma^{-1}$ has the same distribution.
See section 2 in \cite{Pman}, and see
 the final paragraph of section 2 in \cite{P22} where a similar result is proved.)

The motivation behind this note was to find some natural classes of distributions for which the optimal asymptotic probability of success in the secretary problem is less than $\frac1e$. Of course, in light of the discussion in the above paragraph,
the independence noted there must not hold. In this note, for each of the six permutations $\eta\in S_3$,
we consider the secretary problem   with items arriving according to a uniformly random permutation which avoids the pattern $\eta$. The definition of a pattern avoiding permutation is given below.
For such  random permutations, the random variables
$\{I_j\}_{j=1}^n$  defined in the previous paragraph are not  independent, and in certain cases the limiting optimal probability of success is less than $\frac1e$.

Our convention will be that the number $n$ is the highest rank.
For $M\in\{0,1,\cdots, n-1\}$, let $\mathcal{S}(n,M)$ denote the  strategy whereby the observer rejects the first $M$ items, and then selects the first later-arriving item whose rank is higher than that of any
of the first $M$ items (if such an item exists). Let $W_{\mathcal{S}(n,M)}$ denote the event that the observer successfully selects the item of largest rank when using strategy $\mathcal{S}(n,M)$.

Before turning to pattern-avoiding permutations, we begin by noting a cooked-up family of distributions that are the most adversarial as  possible.
Note that in the  uniform case, for any $j\in[n]$, if one uses the absurd strategy of choosing the $j$th item no matter what has occurred, then the probability of
success is $\frac1n$.
There are certain distributions for which
any strategy of the form $\mathcal{S}(n,M)$ also yields this low probability of success.
\begin{proposition}\label{low}
Let $P_n^{\text{low}}$ be any distribution on
$S_n$ that satisfies the following conditions:

\noindent i. $P_n^{\text{low}}(\sigma_j=n)=\frac1n,\ j\in[n]$;

\noindent ii. $P_n^{\text{low}}(\sigma_1<\cdots< \sigma_{j-1}|\sigma_j=n)=1,\ j\in[n]$.

\noindent Then $P_n^{\text{low}}(W_{\mathcal{S}(n,M)})=\frac1n,\ \text{for all}\  M\in[n]\ \text{and}\ n\in\mathbb{N}$.
\end{proposition}
\bf\noindent Remark.\rm\ There are many ways to realize such a probability measure $P_n^{\text{low}}$.
For example, for $j\in[n]$, let $\sigma^{n;j}\in S_n$ be any particular permutation satisfying
$$
\sigma^{n;j}_i=\begin{cases} i,\ i<j;\\ n,\ i=j.\end{cases}
$$
Then one can choose $P_n^{\text{low}}$ to be the uniform probability on $\{\sigma^{n;j}\}_{j=1}^n$.
\begin{proof}Let $M\in\{0,\cdots, n-1\}$.
By the definition of strategy $\mathcal{S}(n,M)$,
$P_n^{\text{low}}(W_{\mathcal{S}(n,M)}|\sigma_j=n)=0$, if $1\le j\le M$. On the other hand, if $M+1\le j\le n$, then from the definition of strategy
$\mathcal{S}(n,M)$ and the second condition on $P_n^{\text{low}}$ in the statement of the proposition,
it follows that
$$
P_n^{\text{low}}(W_{\mathcal{S}(n,M)}|\sigma_j=n)=\begin{cases}1,\ \text{if}\ j=M+1;\\ 0,\ \text{if}\ M+2\le j\le n.\end{cases}
$$
Thus, $P_n^{\text{low}}(W_{\mathcal{S}(n,M)})=P_n^{\text{low}}(\sigma_{M+1}=n)=\frac1n$,
where the final equality follows from the first condition on $P_n^{\text{low}}$ in the statement of the proposition.
\end{proof}

We now turn to pattern-avoiding permutations.
If $\sigma=\sigma_1\cdots\sigma_n\in S_n$ and $\eta=\eta_1\cdots\eta_m\in S_m$, where $2\le m\le n$,
then we say that $\sigma$ contains $\eta$ as a pattern if there exists a subsequence $1\le i_1<i_2<\cdots<i_m\le n$ such
that for all $1\le j,k\le m$, the inequality $\sigma_{i_j}<\sigma_{i_k}$ holds if and only if the inequality $\eta_j<\eta_k$ holds.
If $\sigma$ does not contain $\eta$, then we say that $\sigma$ \it avoids\rm\ $\eta$.
We denote by $S_n^{\text{av}(\eta)}$ the set of permutations in $S_n$ that avoid $\eta$.
%If $n<m$, we define $S_n^{\text{av}(\eta)}=S_n$.

For any $\eta\in S_3$, the number of $\eta$-avoiding permutations is counted by the Catalan numbers: $|S_n^{\text{av}(\eta)}|=C_n$, where $C_n$, the $n$th Catalan number, is given by
$C_n=\frac1{n+1}\binom{2n}n$,\ $n\in\mathbb{N}$ \cite{B04}. We also define $C_0=1$.

For any $\eta\in S_3$, we denote the uniform probability measure on $S_n^{\text{av}(\eta)}$ by $P_n^{\text{av}(\eta)}$.
In all of the results below, we assume tacitly that $n\ge3$.
Our main result concerns the secretary problem with items arriving according to a random permutation distributed according to $P_n^{\text{av}(\eta)}$  in the case that $\eta=231$ or $\eta=132$.

\begin{theorem}\label{1} Let $\eta\in\{231, 132\}$.
For any choice of  $M\in\{0,\cdots, n-1\}$,
\begin{equation}\label{result}
P_n^{\text{av}(\eta)}(W_{\mathcal{S}(n,M)})=\frac{C_{n-1}}{C_n}=\frac{n+1}{2(2n-1)}.
\end{equation}
In particular, for any choice of $M_n\in\{0,\cdots, n-1\}$, the optimal probability
$P_n^{\text{av}(\eta)}(W_{\mathcal{S}(n,M_n)})$ is decreasing in $n$ and
$\lim_{n\to\infty}P_n^{\text{av}(231)}(W_{\mathcal{S}(n,M_n)})=\frac14$.
\end{theorem}
\noindent \bf Remark.\rm\ It is rather anomalous that $P_n^{\text{av}(\eta)}(W_{\mathcal{S}(n,M)})$ is independent of $M\in\{0,\cdots, n-1\}$, for $\eta\in\{231,132\}$.
\medskip

We also have  complete results for the cases $\eta=123$ and $\eta=213$.

\begin{proposition}\label{123}
\begin{equation}\label{123result}
\begin{aligned}
&P_n^{\text{av}(123)}(W_{\mathcal{S}(n,0)})=P_n^{\text{av}(123)}(\sigma_1=n)=\frac{C_{n-1}}{C_n};\\
&P_n^{\text{av}(123)}(W_{\mathcal{S}(n,M)})=P_n^{\text{av}(123)}(n\in\{\sigma_{M+1},\cdots,\sigma_n\}),\ M\in\{1,\cdots, n-1\}.
\end{aligned}
\end{equation}
In particular, the optimal probability occurs for $M=1$ and is increasing in $n$.  One has
$$
P_n^{\text{av}(\eta)}(W_{\mathcal{S}(n,1)})=P_n^{\text{av}(\eta)}(\sigma_1\neq n)=1-\frac{C_{n-1}}{C_n}
$$
and
$$
\lim_{n\to\infty}P_n^{\text{av}(231)}(W_{\mathcal{S}(n,1)})=\frac34.
$$
\end{proposition}

\begin{proposition}\label{213}
For   $M\in\{0,\cdots, n-1\}$,
\begin{equation}\label{213result}
P_n^{\text{av}(213)}(W_{\mathcal{S}(n,M)})=P_n^{\text{av}(213)}(\sigma_{M+1}=n).
\end{equation}
The   probability  $P_n^{\text{av}(213)}(\sigma_{M+1}=n)$ is non-increasing in $M$, and its maximum occurs for $M\in\{0,1\}$. One has
$$
P_n^{\text{av}(213)}(W_{\mathcal{S}(n,0)})=P_n^{\text{av}(213)}(W_{\mathcal{S}(n,1)})=\frac{C_{n-1}}{C_n}=\frac{n+1}{2(2n-1)}.
$$
In particular, the optimal probability is decreasing in $n$ and
$$
\lim_{n\to\infty}P_n^{\text{av}(213)}(W_{\mathcal{S}(n,0)})=\lim_{n\to\infty}P_n^{\text{av}(213)}(W_{\mathcal{S}(n,1)})=\frac14.
$$
\end{proposition}

For the other two permutations in $S_3$, we don't have exact results, but it is easy to show that the limiting optimal  probability of success
is greater than $\frac1e$.
\begin{proposition}\label{othercases}
Let $\eta\in\{321,312\}$. Then
\begin{equation}\label{othercasesresult}
\liminf_{n\to\infty}\max_{M\in\{0,\cdots, n-1\}}P_n^{\text{av}(\eta)}(W_{\mathcal{S}(n,M)})\ge\frac7{16}.
\end{equation}
\end{proposition}

 Theorem \ref{1} is proved in section \ref{secproof}, and Propositions \ref{123}, \ref{213} and \ref{othercases} are proved in section \ref{otherproofs}.

\section{Proof of Theorem \ref{1}}\label{secproof}
As noted in the introduction, for any $\eta\in S_3$, $|S_n^{\text{av}(\eta)}|=C_n$. If $\sigma\in S_n^{\text{av}(231)}$ and $\sigma_j=n$, then
$\{\sigma_1,\cdots,\sigma_{j-1}\}=\{1,\cdots, j-1\}$ and $\{\sigma_{j+1},\cdots\sigma_n\}=\{j,\cdots, n-1\}$.
Similarly, if $\sigma\in S_n^{\text{av}(132)}$ and $\sigma_j=n$, then
$\{\sigma_1,\cdots,\sigma_{j-1}\}=\{n-j+1,\cdots, n-1\}$ and $\{\sigma_{j+1},\cdots\sigma_n\}=\{1,\cdots, n-j\}$.
From this, one obtains the well-known recursion formula for Catalan numbers:
\begin{equation}\label{catrecurr}
C_n=\sum_{j=1}^n C_{j-1}C_{n-j},\ n\in\mathbb{N},
\end{equation}
as well as the formula
\begin{equation}\label{probn}
P_n^{\text{av}(231)}(\sigma_j=n)=P_n^{\text{av}(132)}(\sigma_j=n)=
\frac{C_{j-1}C_{n-j}}{C_n},\ j\in[n].
\end{equation}

We will prove the theorem for the case of 231-avoiding permutations, the proof for 132-avoiding permutations being virtually the same.
Our setup is that $n$ ranked items, labelled from 1 to $n$, arrive in a random permuted order, where the random permutation is distributed according to the probability measure
$P_n^{\text{av}(231)}$ on $S_n$. The highest ranked item corresponds to the number $n$.

Consider first the case that $M=0$. If the observer uses strategy $\mathcal{S}(n,0)$, then the  observer successfully selects the item of highest rank if and only if
$\sigma_1=n$. Thus, from \eqref{probn},
\begin{equation}\label{M=0}
P_n^{\text{av}(231)}(W_{\mathcal{S}(n,0)})=P_n^{\text{av}(231)}(\sigma_1=n)=\frac{C_0C_{n-1}}{C_n}=\frac{n+1}{2(2n-1)}.
\end{equation}

Now consider the case $M\in\{1,\cdots, n-1\}$. If the observer uses strategy $\mathcal{S}(n,M)$, and if
$\sigma_j=n$, then according to the strategy, the observer automatically loses if $j\le M$, while if $j\in\{M+1,\cdots, n\}$, then the observer successfully selects the item of highest rank if and only if
$\max\{\sigma_1,\cdots, \sigma_M\}=\max\{\sigma_1,\cdots, \sigma_{j-1}\}$. Thus,
\begin{equation}\label{basicprobform}
\begin{aligned}
&P_n^{\text{av}(231)}(W_{\mathcal{S}(n,M)})=\sum_{j=M+1}^nP_n^{\text{av}(231)}(\sigma_j=n)P_n^{\text{av}(231)}(W_{\mathcal{S}(n,M)}|\sigma_j=n)=\\
&\sum_{j=M+1}^nP_n^{\text{av}(231)}(\sigma_j=n)P_n^{\text{av}(231)}(\max\{\sigma_1,\cdots, \sigma_M\}=\max\{\sigma_1,\cdots, \sigma_{j-1}\}|\sigma_j=n).
\end{aligned}
\end{equation}
%From \eqref{probn}, we have
As noted at the beginning of the  section, conditioned on $\{\sigma_j=n\}$, we have $\{\sigma_1,\cdots, \sigma_{j-1}\}=\{1,\cdots, j-1\}$.
Thus, under the conditional probability $P_n^{\text{av}(231)}(\ \cdot\ |\sigma_j=n)$, $\sigma_1\cdots\sigma_{j-1}$ is distributed according to
$P_{j-1}^{\text{av}(231)}$. Using this and applying \eqref{probn} with $j-1$ in  place of $n$ and $i$ in place of $j$, we have
\begin{equation}\label{key}
\begin{aligned}
&P_n^{\text{av}(231)}(\max\{\sigma_1,\cdots, \sigma_M\}=\max\{\sigma_1,\cdots, \sigma_{j-1}\}|\sigma_j=n)=\\
&P_{j-1}^{\text{av}(231)}(\max\{\sigma_1,\cdots, \sigma_M\}=\max\{\sigma_1,\cdots, \sigma_{j-1}\})=\\
&\sum_{i=1}^{M}P_{j-1}^{\text{av}(231)}(\max(\sigma_1,\cdots, \sigma_{j-1})=\sigma_i)=\sum_{i=1}^M\frac{C_{i-1}C_{j-1-i}}{C_{j-1}}.
\end{aligned}
\end{equation}
From \eqref{basicprobform}, \eqref{key} and\eqref{probn}, we conclude that
\begin{equation}\label{finalprobform}
P_n^{\text{av}(231)}(W_{\mathcal{S}(n,M)})=\sum_{j=M+1}^n \sum_{i=1}^M\frac{C_{n-j}}{C_n}C_{i-1}C_{j-1-i},\ 1\le M\le n-1.
\end{equation}

For convenience, let
\begin{equation}\label{convenience}
D_{n,M}=\sum_{j=M+1}^n \sum_{i=1}^M\frac{C_{n-j}}{C_n}C_{i-1}C_{j-1-i},\ 1\le M\le n-1.
\end{equation}
For $1\le M\le n-2$, we have
\begin{equation}\label{Dnm}
\begin{aligned}
&D_{n,M+1}-D_{n,M}=\sum_{j=M+2}^n \sum_{i=1}^{M+1}\frac{C_{n-j}}{C_n}C_{i-1}C_{j-1-i}-\sum_{j=M+2}^n \sum_{i=1}^M\frac{C_{n-j}}{C_n}C_{i-1}C_{j-1-i}+\\
&\sum_{j=M+2}^n \sum_{i=1}^M\frac{C_{n-j}}{C_n}C_{i-1}C_{j-1-i}-\sum_{j=M+1}^n \sum_{i=1}^M\frac{C_{n-j}}{C_n}C_{i-1}C_{j-1-i}=\\
&\sum_{j=M+2}^n\frac{C_{n-j}}{C_n}C_MC_{j-2-M}-\sum_{i=1}^M\frac{C_{n-M-1}}{C_n}C_{i-1}C_{M-i}.
\end{aligned}
\end{equation}
From \eqref{catrecurr},
\begin{equation}\label{2catrecurr}
\begin{aligned}
&\sum_{j=M+2}^nC_{n-j}C_{j-2-M}=\sum_{l=1}^{n-M-1}C_{l-1}C_{n-M-1-l}=C_{n-M-1};\\
&\sum_{i=1}^MC_{i-1}C_{M-i}=C_M.
\end{aligned}
\end{equation}
From \eqref{Dnm} and \eqref{2catrecurr}, we have
\begin{equation}\label{zero}
D_{n,M+1}-D_{n,M}=\frac{C_MC_{n-M-1}}{C_n}-\frac{C_{n-M-1}C_M}{C_n}=0,\ 1\le M\le n-2.
\end{equation}
We also have
\begin{equation}\label{Dn1}
D_{n,1}=\sum_{j=2}^n\frac{C_{n-j}}{C_n}C_0C_{j-2}=\frac1{C_n}\sum_{l=1}^{n-1}C_{n-1-l}C_{l-1}=\frac{C_{n-1}}{C_n}=\frac{n+1}{2(2n-1)}.
\end{equation}
Now \eqref{result} follows from \eqref{M=0}, \eqref{finalprobform}, \eqref{convenience}, \eqref{zero} and \eqref{Dn1}.
The final statement in the theorem follows immediately from \eqref{result}.
\hfill $\square$

\section{Proofs of Propositions \ref{123}, \ref{213} and \ref{othercases}}\label{otherproofs}
\noindent \it Proof of Proposition \ref{123}.\rm\
If $\sigma\in S_n$ satisfies $\sigma_1=n$, then $\sigma\in S_n^{\text{av}(123)}$ if and only if $\sigma_2\cdots\sigma_n\in S_{n-1}^{\text{av}(123)}$. Thus,
$P_n^{\text{av}(123)}(\sigma_1=n)=\frac{C_{n-1}}{C_n}$. Consequently, by the definition of strategy $\mathcal{S}(n,0)$,
$P(W_{\mathcal{S}(n,0)})=P_n^{\text{av}(123)}(\sigma_1=n)=\frac{C_{n-1}}{C_n}$.

If $\sigma\in S_n^{\text{av}(123)}$ and
$\sigma_j=n$, for $j\in \{2,\cdots,n\}$, then $\sigma_1>\cdots>\sigma_{j-1}$.
Consequently, for $M\in\{1,\cdots, n-1\}$,  by the definition of the strategy $\mathcal{S}(n,M)$,
 $P_n^{\text{av}(123)}(W_{\mathcal{S}(n,M)})=P_n^{\text{av}(123)}(n\in\{\sigma_{M+1},\cdots,\sigma_n\})$. From this and the fact that $\frac{C_{n-1}}{C_n}=\frac{n+1}{2(2n-1)}<\frac12$ (as noted before the statement of the results, we are assuming $n\ge3$), it follows that the optimal strategy is
 $\mathcal{S}(n,1)$ and that  $P_n^{\text{av}(123)}(W_{\mathcal{S}(n,1)})=P(\sigma_1\neq n)=1-\frac{C_{n-1}}{C_n}$.
The rest of the proposition follows from the fact that $\frac{C_{n-1}}{C_n}=\frac{n+1}{2(2n-1)}$.
\hfill $\square$
\medskip

\noindent \it Proof of Proposition \ref{213}.\rm\
Let $M\in\{0,\cdots, n-1\}$.
By the definition of strategy $\mathcal{S}(n,M)$,
$P_n^{\text{av}(213)}(W_{\mathcal{S}(n,M)}|\sigma_j=n)=0$, if $1\le j\le M$.
Now consider the case that $\sigma_j=n$ with $M+1\le j\le n$. By the definition of 213-pattern avoiding,
$$
P_n^{\text{av}(213)}(\sigma_1<\sigma_2<\cdots<\sigma_{j-1}|\sigma_j=n)=1.
$$
Using this with the definition of the strategy $\mathcal{S}(n,M)$, it follows that
$$
P_n^{\text{av}(213)}(W_{\mathcal{S}(n,M)}|\sigma_j=n)=\begin{cases} 1,\ \text{if}\ j=M+1;\\ 0,\ \text{if}\ M+2\le j\le n.\end{cases}
$$
From the above, we conclude that
$$
P_n^{\text{av}(213)}(W_{\mathcal{S}(n,M)})=P_n^{\text{av}(213)}(\sigma_{M+1}=n).
$$

For any $j\in\{2,\cdots, n\}$ and  any permutation $\sigma\in S_n^{\text{av}(213)}$ which satisfies $\sigma_j=n$,
it is easy to see that  $\sigma'\in S_n^{\text{av}(213)}$, where
$$
\sigma'_i=\begin{cases}\sigma_i,\  i\not\in \{j-1,j\};\\ n,\ i=j-1;\\ \sigma_{j-1},\ i=j.\end{cases}
$$
This shows that $P_n^{\text{av}(213)}(\sigma_{M+1}=n)$ is non-increasing in $M$.
It is easy to see that for any permutation $\sigma\in S_n$ satisfying $\sigma_1=n$,
one has $\sigma\in S_n^{\text{av}(213)}$ if and only if
the permutation $\sigma_2\cdots\sigma_n\in S_{n-1}^{\text{av}(213)}$,
and that for any $\sigma\in S_n$ satisfying $\sigma_2=n$,
one has $\sigma\in S_n^{\text{av}(213)}$ if and only if
the permutation $\sigma_1\sigma_3\cdots\sigma_n\in S_{n-1}^{\text{av}(213)}$.
 Thus,
$P_n^{\text{av}(213)}(\sigma_1=n)=P_n^{\text{av}(213)}(\sigma_2=n)=
\frac{C_{n-1}}{C_n}$.
\hfill $\square$
\medskip

\noindent \it Proof of Proposition \ref{othercases}.\rm\
We consider the cases $\eta=321$ and $\eta=312$ together.
By the definition of the strategy $\mathcal{S}(n,n-2)$, the event $W_{\mathcal{S}(n,n-2)}$ of successfully selecting the highest ranked item when using strategy  $\mathcal{S}(n,n-2)$ satisfies
\begin{equation}\label{eventn-2}
\{\sigma_{n-1}=n\}\cup\{\sigma_{n-1}\neq n-1,\sigma_n=n\}= W_{\mathcal{S}(n,n-2)}.
\end{equation}
For $\eta=321$ or $\eta=312$, a permutation $\sigma\in S_n$
 satisfying $\sigma_{n-1}=n$ belongs to $S_n^{\text{av}(\eta)}$ if and only if the permutation $\sigma_1\cdots\sigma_{n-2}\sigma_n$ belongs to $S_{n-1}^{\text{av}(\eta)}$.
Similarly, a permutation $\sigma\in S_n$
 satisfying $\sigma_n=n$ belongs to $S_n^{\text{av}(\eta)}$ if and only if the permutation $\sigma_1\cdots\sigma_{n-1}$ belongs to $S_{n-1}^{\text{av}(\eta)}$.
Also, a permutation $\sigma\in S_n$
 satisfying $\sigma_{n-1}=n-1$ and $\sigma_n=n$ belongs to $S_n^{\text{av}(\eta)}$ if and only if the permutation $\sigma_1\cdots\sigma_{n-2}$ belongs to $S_{n-2}^{\text{av}(\eta)}$.
From this it follows that
\begin{equation}\label{calc}
\begin{aligned}
&P_n^{\text{av}(\eta)}(\sigma_n=n)=
P_n^{\text{av}(\eta)}(\sigma_{n-1}=n)=
\frac{C_{n-1}}{C_n};\\
&P_n^{\text{av}(\eta)}(\sigma_{n-1}=n-1,\sigma_n=n)=\frac{C_{n-2}}{C_n},\\
&\text{for}\ \eta\in\{321,312\}.
\end{aligned}
\end{equation}
From \eqref{eventn-2} and\eqref{calc} it follows that
$$
P_n^{\text{av}(\eta)}(W_{\mathcal{S}(n,n-2)})=2\frac{C_{n-1}}{C_n}-\frac{C_{n-2}}{C_n}, \ \text{for}\ \eta\in\{321,312\}.
$$
Using this with the well-know asymptotic formula, $C_n\sim \frac{4^n}{\sqrt\pi\thinspace n^\frac32}$,
we have
$$
\lim_{n\to\infty}P_n^{\text{av}(\eta)}(W_{\mathcal{S}(n,n-2)})= \frac7{16},\ \text{for}\ \eta\in\{321,312\},
$$
which proves the proposition.
\hfill $\square$

\end{document}